# Analyzing Incomplete Discrete Longitudinal Clinical Trial Data

Ivy Jansen, Caroline Beunckens, Geert Molenberghs, Geert Verbeke and Craig Mallinckrodt

*Abstract.* Commonly used methods to analyze incomplete longitudinal clinical trial data include complete case analysis (CC) and last observation carried forward (LOCF). However, such methods rest on strong assumptions, including missing completely at random (MCAR) for CC and unchanging profile after dropout for LOCF. Such assumptions are too strong to generally hold. Over the last decades, a number of full longitudinal data analysis methods have become available, such as the linear mixed model for Gaussian outcomes, that are valid under the much weaker missing at random (MAR) assumption. Such a method is useful, even if the scientific question is in terms of a single time point, for example, the last planned measurement occasion, and it is generally consistent with the intention-to-treat principle. The validity of such a method rests on the use of maximum likelihood, under which the missing data mechanism is ignorable as soon as it is MAR. In this paper, we will focus on non-Gaussian outcomes, such as binary, categorical or count data. This setting is less straightforward since there is no unambiguous counterpart to the linear mixed model. We first provide an overview of the various modeling frameworks for non-Gaussian longitudinal data, and subsequently focus on generalized linear mixed-effects models, on the one hand, of which the parameters can be estimated using full likelihood, and on generalized estimating equations, on the other hand, which is a nonlikelihood method and hence requires a modification to be valid under MAR. We briefly comment on the position of models that assume missingness not at random and argue they are most useful to perform sensitivity analysis. Our developments are underscored using data from two studies. While the case studies feature binary outcomes, the methodology applies equally well to other discrete-data settings, hence the qualifier "discrete" in the title.

Ivy Jansen is Postdoctoral Researcher, Center for Statistics, Hasselt University, Agoralaan Gebouw D, B-3590 Diepenbeek, Belgium (e-mail: ivy.jansen@uhasselt.be). Caroline Beunckens is Research Assistant, Center for Statistics, Hasselt University, Agoralaan Gebouw D, B-3590 Diepenbeek, Belgium (e-mail: caroline.beunckens@uhasselt.be). Geert Molenberghs is Professor, Center for Statistics, Hasselt University, Agoralaan Gebouw D, B-3590 Diepenbeek, Belgium (e-mail: geert.molenberghs@uhasselt.be).

## 1. INTRODUCTION

Data from longitudinal studies, in general, and from clinical trials, in particular, are prone to incompleteness. Dropout is a special case of incompleteness. Since incompleteness usually occurs for reasons outside the control of the investigators and may be related to the outcome measurement of interest, it is generally necessary to address the process that governs incompleteness. Only in special but important cases is it possible to ignore the missingness process.

When referring to the missing-value, or nonresponse, process, we will use the terminology of Little and Rubin (2002, Chapter 6). A nonresponse process is said to be *missing completely at random* (MCAR) if the missingness is independent of both unobserved and observed data, and said to be *missing at random* (MAR) if, conditional on the observed data, the missingness is independent of the unobserved measurements. A process that is neither MCAR nor MAR is termed *nonrandom* (MNAR). In the context of likelihood inference, and when the parameters that describe the measurement process are functionally independent of the parameters that describe the missingness process, MCAR and MAR are *ignorable*, while a nonrandom process is nonignorable.

Early work regarding missingness focused on the consequences of the induced lack of balance of deviations from the study design (Afifi and Elashoff, 1966; Hartley and Hocking, 1971). Later, algorithmic developments took place, such as the expectation–maximization algorithm (EM; Dempster, Laird and Rubin, 1977) and multiple imputation (Rubin, 1987). These advances have brought likelihood-based ignorable analysis within reach for a large class of designs and models. However, they usually require extra programming in addition to available standard statistical software.

*Geert Verbeke is Professor, Biostatistical Centre, Katholieke Universiteit Leuven, Kapucijnenvoer 35, B-3000 Leuven, Belgium (e-mail: geert.verbeke@med.kuleuven.be). Craig Mallinckrodt is Senior Research Fellow, Eli Lilly & Company, Indianapolis, Indiana 46285, USA (e-mail: mallinckrodt_craig@Lilly.com).*

In the meantime, however, clinical trial practice has put a strong emphasis on such methods as *complete case analysis* (CC), which restricts the analysis to those subjects for which all information has been measured according to protocol, and *last observation carried forward* (LOCF), for which the last observed measurement is substituted for values at later points in time that are not observed, or other simple forms of imputation. Claimed advantages include computational simplicity, no need for a full longitudinal model (e.g., when the scientific question is in terms of the last planned measurement occasion only) and, for LOCF, compatibility with the intention-to-treat (ITT) principle. Within the Gaussian setting, Molenberghs et al. (2004) have argued that this focus is understandable, but, given current computational resources, unfortunate. They suggest the use of a likelihood-based ignorable analysis, for example, based on the linear mixed-effects model. Such a method requires MAR rather than the much stronger assumptions that underlie CC and LOCF, and uses all data, obviating the need for both deleting and filling in data, and is thus consistent with the intention-to-treat principle. Nevertheless, care has to be taken when subjects are discontinued for reasons of noncompliance, since then the modes of analysis indicated here would assume that treatment is unchanged after dropout. This implies the need for sensitivity analysis and an important discussion of this point has been given by Fitzmaurice (2003). Molenberghs et al. (2004) also show that the incomplete sequences contribute to estimands of interest, even early dropouts when scientific interest is in the last planned measurement only. Finally, they show that such an analysis is possible, without the need for any additional data manipulation, using, for example, the SAS procedure MIXED or the SPlus or R function lme. Of course, a longitudinal model has to be specified for the entire vector of responses. In a clinical trial setting, with relatively short and balanced response sequences, full multivariate models, encompassing full treatment-by-group interactions, perhaps corrected for baseline covariates, and an unstructured variance–covariance matrix, are usually within reach. A model of this type is relatively mild in the restrictions made.

The non-Gaussian setting is different in the sense that there is no generally accepted counterpart to the linear mixed-effects model. We therefore first sketch a general taxonomy for longitudinal models in this context, including marginal, random-effects (or subject-specific) and conditional models. We then argue that marginal and random-effects models both have their merits in the analysis of longitudinal clinical trial data and we focus on two important representatives, that is, the generalized estimating equations (GEE) approach within the marginal family and the generalized linear mixed-effects model (GLMM) within the random-effects family. We highlight important similarities and differences between these model families. While GLMM parameters can be fitted using maximum likelihood, the same is not true for the GEE method, which is of a frequentist nature. Therefore, Robins, Rotnitzky and Zhao (1995) have devised so-called weighted generalized estimating equations (WGEE), which are valid under MAR but require the specification of a dropout model in terms of observed outcomes and/or covariates, in view of specifying the weights. Thus, we believe that, generally, methods such as complete case analysis or LOCF ought to be abandoned in favor of the likelihood-based and weighted GEE models discussed here.

By definition, MNAR missingness cannot be fully ruled out based on the observed data. Nevertheless, ignorable analyses may provide reasonably stable results, even when the assumption of MAR is violated, in the sense that such analyses constrain the behavior of the unseen data to be similar to that of the observed data. A discussion of this phenomenon in the survey context has been given in Rubin, Stern and Vehovar (1995). These authors first argue that, in well-conducted experiments (some surveys and many confirmatory clinical trials), the assumption of MAR is often to be regarded as a realistic one. Second, and very important for confirmatory trials, an MAR analysis can be specified a priori without additional work relative to a situation with complete data. Third, while MNAR models are more general and explicitly incorporate the dropout mechanism, the inferences they produce are typically highly dependent on the untestable and often implicit built-in assumptions regarding the distribution of the unobserved measurements given the observed ones. The quality of the fit to the observed data need not reflect at all the appropriateness of the implied structure that governs the unobserved data. Based on these considerations, we recommend, for primary analysis purposes, the use of ignorable likelihood-based methods or appropriately modified frequentist methods. To explore the impact of deviations from the MAR assumption on the conclusions, one should ideally conduct a sensitivity analysis (Verbeke and Molenberghs, 2000, Chapters 18–20).

Two case studies motivate our work. The first one arises from a randomized, double-blind psychiatric clinical trial conducted in the United States. The primary objective of this trial was to compare the efficacy of an experimental anti-depressant with placebo in order to support a new drug application. The study enrolled 167 patients. The Hamilton depression rating scale ($HAMD_{17}$) is used to measure the depression status of the patients. The binary indicator of interest is 1 if the $HAMD_{17}$ score is greater than 7, and 0 otherwise. For each patient, a baseline assessment is available, as well as eight post-baseline visits going from visit 4 to visit 11.

The second case study arises from a randomized multicentric clinical trial that compared an experimental treatment (interferon-$\alpha$) to a corresponding placebo in the treatment of patients with age-related macular degeneration. Interest focuses on the comparison between placebo and the highest dose (6 million units daily) of interferon-$\alpha$ ($Z$), but the full results of this trial have been reported elsewhere (Pharmacological Therapy for Macular Degeneration Study Group, 1997). Patients with macular degeneration progressively lose vision. In the trial, the patients' visual acuity was assessed at different time points (4 weeks, 12 weeks, 24 weeks and 52 weeks) through their ability to read lines of letters on standardized vision charts. These charts display lines of five letters of decreasing size, which the patient must read from top (largest letters) to bottom (smallest letters). Each line with at least four letters correctly read is called one line of vision. The patient's visual acuity is the total number of letters correctly read. The primary endpoint of the trial was the loss of at least three lines of vision at 1 year, compared to their baseline performance (i.e., a binary endpoint). The total number of longitudinal profiles is 240, but only for 188 of these have the four follow-up measurements been made. An overview is given in Table 1. Thus, 78.33% of the profiles are complete, while 18.33% exhibit monotone missingness. Out of the latter group, 2.5% or six subjects have no follow-up measurements. The remaining 3.33%, representing eight subjects, have intermittent missing values. Although the group of dropouts is of considerable magnitude, the ones with intermittent missingness is much smaller. Nevertheless, it is cautious to include all data in the analyses. Data on this second case study are available on the author's website.

TABLE 1
*Age related macular degeneration trial: Overview of missingness patterns and the frequencies with which they occur* (O *indicates observed and* M *indicates missing*)

| **Measurement occasion (weeks)** | | | | | |
|---|---|---|---|---|---|
| **4** | **12** | **24** | **52** | **Number** | **%** |
| | | | Completers | | |
| O | O | O | O | 188 | 78.33 |
| | | | Dropouts | | |
| O | O | O | M | 24 | 10.00 |
| O | O | M | M | 8 | 3.33 |
| O | M | M | M | 6 | 2.50 |
| M | M | M | M | 6 | 2.50 |
| | | Nonmonotone missingness | | | |
| O | O | M | O | 4 | 1.67 |
| O | M | M | O | 1 | 0.42 |
| M | O | O | O | 2 | 0.83 |
| M | O | M | M | 1 | 0.42 |

The general data setting is introduced in Section 2, as well as a formal framework for incomplete longitudinal data, together with a discussion of the problems associated with simple methods. Section 3 focuses on two important families of models for discrete repeated measures. The first case study is analyzed in Section 4, while the second one is the subject of Section 5. The paper ends with a discussion in Section 6.

## 2. DATA SETTING AND MODELING FRAMEWORK

Assume that for subject $i = 1, \ldots, N$ a sequence of responses $Y_{ij}$ is designed to be measured at occasions $j = 1, \ldots, n$. The outcomes are grouped into a vector $\mathbf{Y}_i = (Y_{i1}, \ldots, Y_{in})'$. In addition, for each occasion $j$, define $R_{ij}$ as being equal to 1 if $Y_{ij}$ is observed and 0 otherwise. The missing data indicators $R_{ij}$ are grouped into a vector $\mathbf{R}_i$, which is of the same length as $\mathbf{Y}_i$. Define now a dropout indicator $D_i$ for the occasion at which dropout occurs and make the convention that $D_i = n + 1$ for a complete sequence. Further, split the vector $\mathbf{Y}_i$ into observed ($\mathbf{Y}_i^o$) and missing ($\mathbf{Y}_i^m$) components, respectively.

Modeling usually is initiated by considering the full data density $f(\mathbf{y}_i, d_i | \boldsymbol{\theta}, \boldsymbol{\psi})$, where the parameter vectors $\boldsymbol{\theta}$ and $\boldsymbol{\psi}$ describe the measurement and missingness processes, respectively. Covariates are assumed to be measured, but for notational simplicity are suppressed from notation.

Most strategies used to analyze such data are, implicitly or explicitly, based on the following two choices.

*Model for measurements.* A choice has to be made regarding the modeling approach to the measurements. There are three common views. In the first view, one opts to analyze the entire longitudinal profile, irrespective of whether interest focuses on the entire profile, on the one hand (e.g., difference in slope between groups), or whether a specific time is of interest, on the other hand (e.g., the last planned occasion). In the latter case, the motivation to model the entire profile is because, for example, earlier responses do provide statistical information on later ones. This is especially true when dropout is present. In the second view, one defines the scientific question and restricts the corresponding analysis to the *last planned occasion*. Of course, as soon as dropout occurs, such a measurement may not be available, whence often *last observation carried forward* is used. Such an analysis is based on strong assumptions that often do not hold. This point has been made extensively in Molenberghs et al. (2004). These authors advocate the use of proper, likelihood-based, longitudinal methods, which are generally valid under MAR. In the third view, one chooses to define the question and the corresponding analysis in terms of the *last observed measurement*. While sometimes used as an alternative motivation for so-called *last observation carried forward* analyses (Siddiqui and Ali, 1998; Mallinckrodt, Clark, Carroll and Molenberghs, 2003a; Mallinckrodt et al., 2003b), a common criticism is that the last observed measurement amalgamates measurements at real stopping times (for dropouts) and at a purely design-based time (for completers). Thus, we hope to show that the first view is in many situations the most sensible route of analysis.

*Method for handling missingness.* A choice has to be made regarding the modeling approach to the missingness model. Luckily, under certain assumptions this process can be ignored (likelihood-based or Bayesian ignorable analysis, for which MAR is a sufficient condition). Some simple methods such as CC *analysis* and LOCF do not explicitly address the missingness mechanism either, but are nevertheless not ignorable. We will return to this issue in the next section.

Let us first describe the measurement and missingness models in turn, and then introduce and comment on ignorability. The measurement model will depend on whether or not a full longitudinal analysis is done. In case the second or third view is adopted, one typically opts for classical two-group or multi-group comparisons ($t$ test, Wilcoxon, etc.). In case a longitudinal analysis is deemed necessary, the choice made depends on the nature of the outcome.

For continuous outcomes, a common choice is the general linear mixed-effects model or a special case of it, such as a (structured) multivariate normal model (Verbeke and Molenberghs, 2000). However, for categorical (nominal, ordinal and binary) and discrete outcomes (counts), as in our case study, the modeling choices are less straightforward. Extensions of the generalized linear models to the longitudinal case were discussed by Diggle, Heagerty, Liang and Zeger (2002), where a lot of emphasis is on generalized estimating equations (Liang and Zeger, 1986). Generalized linear mixed models have been proposed and/or studied by, for example, Stiratelli, Laird and Ware (1984), Wolfinger and O'Connell (1993) and Breslow and Clayton (1993). Fahrmeir and Tutz (2001) devoted an entire book to generalized linear models for multivariate settings. We return to modeling non-Gaussian repeated measures in Section 3. It is important to note that, since quite distinct modeling families are in use, the researcher ought to be guided by the main scientific question at hand when choosing between the modeling families.

Assume that incompleteness is due to dropout only and that the first measurement $Y_{i1}$ is obtained for everyone. The model for the dropout process is based on, for example, a logistic regression for the probability of dropout at occasion $j$, given the subject is still in the study. We denote this probability by $g(\mathbf{h}_{ij}, y_{ij})$, in which $\mathbf{h}_{ij}$ is a vector that contains all responses observed up to but not including occasion $j$, as well as relevant covariates. We then assume that $g(\mathbf{h}_{ij}, y_{ij})$ satisfies

$$\begin{aligned}
&\text{logit}[g(\mathbf{h}_{ij}, y_{ij})]\\
(1)\quad &= \text{logit}[\text{pr}(D_i = j | D_i \geq j, \mathbf{y}_i)]\\
&= \mathbf{h}_{ij}\boldsymbol{\psi} + \omega y_{ij}, \qquad i = 1, \ldots, N.
\end{aligned}$$

When $\omega$ equals zero, and assuming the posited model is correct, the dropout model is MAR. If $\omega \neq 0$, the posited dropout process is MNAR. Model (1) provides the building blocks for the dropout process $f(d_i | \mathbf{y}_i, \boldsymbol{\psi})$.

Rubin (1976) and Little and Rubin (2002) have shown that, under MAR and mild regularity conditions (parameters $\boldsymbol{\theta}$ and $\boldsymbol{\psi}$ are functionally independent), likelihood-based and Bayesian inference are valid when the missing data mechanism is ignored (see also Verbeke and Molenberghs, 2000). Practically speaking, the likelihood of interest is then based on the factor $f(\mathbf{y}_i^o | \boldsymbol{\theta})$. This is called *ignorability*. A model of the form (1), of course with $\omega = 0$, may but does not have to be considered in such a case. The practical implication is that a software module with likelihood estimation facilities and with the ability to deal with incompletely observed subjects manipulates the correct likelihood, providing valid parameter estimates, standard errors if based on the observed information matrix and likelihood ratio values (Kenward and Molenberghs, 1998). Examples of such software tools include the MIXED, NLMIXED and GENMOD procedures in SAS.

A few cautionary remarks are in order. First, when at least part of the scientific interest is directed toward the nonresponse process, obviously both processes need to be considered. Still, under MAR, both processes can be modeled and parameters can be estimated separately. Second, it may be hard to fully rule out the operation of an MNAR mechanism. Third, one is now restricted to the first view on modeling the outcomes, that is, a full longitudinal analysis is necessary, even when interest is restricted to, for example, a comparison between the two treatment groups at the last occasion. In the latter case, the fitted model can be used as the basis for inference at the last occasion. A common criticism, especially in a regulated controlled clinical trial setting, is *that a model needs to be considered*. However, it should be noted that in many clinical trial settings the repeated measures are balanced in the sense that a common (and often limited) set of measurement times is considered for all subjects, allowing the a priori (protocol) specification of a saturated model (e.g., full group-by-time interaction model for the means and unstructured variance–covariance matrix).

Such an ignorable linear mixed model specification is termed mixed-effects model repeated measures (MMRM) by Mallinckrodt, Clark and David (2001a, b). Thus, MMRM is a particular form of a linear mixed model, relevant for acute phase confirmatory clinical trials, fitting within the ignorable likelihood paradigm. It has to be noted that this approach, for the special case where no dropout occurs, is fully equivalent to a one-way multivariate analysis of variance model for the repeated outcomes, with a class variable treatment effect. This observation provides a strong basis for such an approach, which is a very promising alternative for the simple ad hoc methods such as complete-case analysis or LOCF. While the above reasoning is tied to the continuous-outcome setting, similar modeling strategies exist for the non-Gaussian case, as discussed in Section 3.

These arguments, supplemented with the availability of software tools within which such multivariate models can be fitted to incomplete data, cast doubts regarding the usefulness of such simple methods as CC and LOCF. This issue has been discussed in detail, in the context of Gaussian outcomes, by Molenberghs et al. (2004). Apart from biases as soon as the missing data mechanism is not MCAR, CC can suffer from severe efficiency losses. Especially since tools have become available to include incomplete sequences along with complete ones into the analysis, one should do everything possible to avoid wasting patient data.

Last observation carried forward, as other imputation strategies (Dempster and Rubin, 1983; Little and Rubin, 2002), can lead to artificially inflated precision. Furthermore as Molenberghs et al. (2004) have shown, the method can produce severely biased treatment comparisons and, perhaps contrary to some common belief, such biases can be conservative but also liberal. The method rests on the strong assumption that a patient's outcome profile remains flat, at the level of the last observed measurement, throughout the remainder of follow-up.

## 3. DISCRETE REPEATED MEASURES

Whereas the linear mixed model and its special cases is seen as a unifying parametric framework for Gaussian repeated measures (Verbeke and Molenberghs, 2000), there are many more options available in the non-Gaussian setting. In a *marginal model*, marginal distributions are used to describe the outcome vector $\mathbf{Y}$, given a set $\mathbf{X}$ of predictor variables. The correlation among the components of $\mathbf{Y}$ can then be captured either by adopting a fully parametric approach or by means of working assumptions, such as in the semiparametric approach of Liang and Zeger (1986). Alternatively, in a *random-effects model*, the predictor variables $\mathbf{X}$ are supplemented with a vector $\boldsymbol{\theta}$ of random effects, conditional upon which the components of $\mathbf{Y}$ are usually assumed to be independent. This does not preclude that more elaborate models are possible if residual dependence is detected (Longford, 1993). Finally, a *conditional model* describes the distribution of the components of $\mathbf{Y}$, conditional on $\mathbf{X}$ but also conditional on (a subset of) the other components of $\mathbf{Y}$. Well-known members of this class of models are log-linear models (Gilula and Haberman, 1994).

Let us give a simple example of each for the case of Gaussian outcomes. A marginal model starts by specifying

$$(2) \qquad E(Y_{ij}|\mathbf{x}_{ij}) = \mathbf{x}'_{ij}\boldsymbol{\beta},$$

whereas in a random-effects model we focus on the expectation, conditional upon the random-effects vector,

$$(3) \qquad E(Y_{ij}|\mathbf{b}_i, \mathbf{x}_{ij}) = \mathbf{x}'_{ij}\boldsymbol{\beta} + \mathbf{z}'_{ij}\mathbf{b}_i.$$

The conditional model uses expectations of the form

$$(4) \quad E(Y_{ij}|Y_{i,j-1},\ldots,Y_{i1},\mathbf{x}_{ij}) = \mathbf{x}'_{ij}\boldsymbol{\beta} + \alpha Y_{i,j-1}.$$

In the linear mixed model case, random-effects models imply a simple marginal model. This is due to the elegant properties of the multivariate normal distribution. In particular, the expectation (2) follows from (3) either by (a) marginalizing over the random effects or by (b) conditioning on the random-effects vector $\mathbf{b}_i = \mathbf{0}$. Hence, the fixed-effects parameters $\boldsymbol{\beta}$ have both a marginal as well as a hierarchical model interpretation.

Since marginal and random-effects models are the most useful ones in our context and given this connection between them, it is clear why the linear mixed model provides a unified framework in the Gaussian setting. Such a close connection between the model families does not exist when outcomes are of a nonnormal type, such as binary, categorical or discrete. We will consider the marginal and random-effects model families in turn and then point to some particular issues that arise within them or when comparisons are made between them. The conditional models are less useful in the context of longitudinal data and will not be discussed here (Molenberghs and Verbeke, 2005).

### 3.1 Marginal Models

Thorough discussions on marginal modeling can be found in Diggle, Heagerty, Liang and Zeger (2002) and Fahrmeir and Tutz (2001). The specific context of clustered binary data has received treatment in Aerts, Geys, Molenberghs and Ryan (2002). Apart from full likelihood approaches, nonlikelihood approaches such as *generalized estimating equations* (Liang and Zeger, 1986) or *pseudolikelihood* (le Cessie and van Houwelingen, 1994; Geys, Molenberghs and Lipsitz, 1998) have been considered.

Bahadur (1961) proposed a marginal model, accounting for the association via marginal correlations. Ekholm (1991) proposed a so-called success probabilities approach. Bowman and George (1995) proposed a model for the particular case of exchangeable binary data. Ashford and Sowden (1970) considered the multivariate probit model for repeated ordinal data, thereby extending univariate probit regression. Molenberghs and Lesaffre (1994) and Lang and Agresti (1994) proposed models that parameterize the association in terms of marginal odds ratios. Dale (1986) defined the bivariate global odds ratio model, based on a bivariate Plackett distribution (Plackett, 1965). Molenberghs and Lesaffre (1994, 1999), Lang and Agresti (1994) and Glonek and McCullagh (1995) extended this model to multivariate ordinal outcomes. They generalized the bivariate Plackett distribution in order to establish the multivariate cell probabilities.

While full likelihood methods are appealing because of their flexible ignorability properties (Section 2), their use for non-Gaussian outcomes can be problematic due to prohibitive computational requirements. Therefore, GEE is a viable alternative within this family. Since GEE is motivated by frequentist considerations, the missing data mechanism needs to be MCAR for it to be ignorable. This motivates the proposal of so-called *weighted generalized estimating equations*. We will discuss these in turn.

3.1.1 *Generalized estimating equations.* Generalized estimating equations, useful to circumvent the computational complexity of full likelihood, can be considered whenever interest is restricted to the mean parameters (treatment difference, time evolutions, effect of baseline covariates, etc.). It is rooted in the quasi-likelihood ideas expressed by McCullagh and Nelder (1989). Modeling is restricted to the correct specification of the marginal mean function, together with so-called *working assumptions* about the correlation structure of the vector of repeated measures.

Let us now introduce the classical form of GEE. Note that the score equations to be solved when computing maximum likelihood estimates under a marginal normal model $\mathbf{y}_i \sim N(X_i\boldsymbol{\beta}, V_i)$ are given by

$$(5) \quad \sum_{i=1}^{N} X_i'(A_i^{1/2} C_i A_i^{1/2})^{-1}(\mathbf{y}_i - X_i\boldsymbol{\beta}) = \mathbf{0},$$

in which the marginal covariance matrix $V_i$ has been decomposed in the form $A_i^{1/2} C_i A_i^{1/2}$, where $A_i$ is the matrix with the marginal variances on the main diagonal and zeros elsewhere, and $C_i$ is equal to the marginal correlation matrix. Switching to the non-Gaussian case, the score equations become

$$(6) \quad S(\boldsymbol{\beta}) = \sum_{i=1}^{N} \frac{\partial \boldsymbol{\mu}_i}{\partial \boldsymbol{\beta}'} (A_i^{1/2} C_i A_i^{1/2})^{-1}(\mathbf{y}_i - \boldsymbol{\mu}_i) = \mathbf{0},$$

which are less linear than (5) due to the presence of a link function (e.g., the logit link for binary data) and the mean–variance relationship. Typically the correlation matrix $C_i$ contains a vector $\boldsymbol{\alpha}$ of unknown parameters that is replaced for practical purposes by a consistent estimate.

Assuming that the marginal mean $\boldsymbol{\mu}_i$ has been correctly specified as $h(\boldsymbol{\mu}_i) = X_i\boldsymbol{\beta}$, it can be shown that, under mild regularity conditions, the estimator $\widehat{\boldsymbol{\beta}}$ obtained by solving (6) is asymptotically normally distributed with mean $\boldsymbol{\beta}$ and with covariance matrix

$$(7) \quad I_0^{-1} I_1 I_0^{-1},$$

where

$$I_0 = \left(\sum_{i=1}^{N} \frac{\partial \boldsymbol{\mu}_i'}{\partial \boldsymbol{\beta}} V_i^{-1} \frac{\partial \boldsymbol{\mu}_i}{\partial \boldsymbol{\beta}'}\right),$$

$$I_1 = \left(\sum_{i=1}^{N} \frac{\partial \boldsymbol{\mu}_i'}{\partial \boldsymbol{\beta}} V_i^{-1} \operatorname{Var}(\mathbf{y}_i) V_i^{-1} \frac{\partial \boldsymbol{\mu}_i}{\partial \boldsymbol{\beta}'}\right).$$

In practice, $\operatorname{Var}(\mathbf{y}_i)$ in (7) is replaced by $(\mathbf{y}_i - \boldsymbol{\mu}_i) \cdot (\mathbf{y}_i - \boldsymbol{\mu}_i)'$, which is unbiased on the sole condition of correct mean specification. One also needs estimates of the nuisance parameters $\boldsymbol{\alpha}$. Liang and Zeger (1986) proposed moment-based estimates for the working correlation. To this end, define deviations

$$e_{ij} = \frac{y_{ij} - \mu_{ij}}{\sqrt{v(\mu_{ij})}}.$$

Some of the more popular choices for the working correlations are independence [$\operatorname{Corr}(Y_{ij}, Y_{ik}) = 0, j \neq k$], exchangeability [$\operatorname{Corr}(Y_{ij}, Y_{ik}) = \alpha, j \neq k$], AR(1) [$\operatorname{Corr}(Y_{ij}, Y_{i,j+t}) = \alpha^t, t = 0, 1, \ldots, n_i - j$] and unstructured [$\operatorname{Corr}(Y_{ij}, Y_{ik}) = \alpha_{jk}, j \neq k$]. Typically, moment-based estimation methods are used to estimate these parameters as part of an integrated iterative estimation procedure. An overdispersion parameter could be included as well, but we have suppressed it for ease of exposition. The standard iterative procedure to fit GEE, based on Liang and

8 JANSEN ET AL.

Zeger (1986), is then as follows: (1) compute initial estimates for $\boldsymbol{\beta}$ using a univariate generalized linear model (i.e., assuming independence); (2) compute the quantities $\mathbf{b}_i$ needed in the estimating equation; (3) compute Pearson residuals $e_{ij}$; (4) compute estimates for $\boldsymbol{\alpha}$; (5) compute $C_i(\boldsymbol{\alpha})$; (6) compute $V_i(\boldsymbol{\beta}, \boldsymbol{\alpha}) = A_i^{1/2}(\boldsymbol{\beta}) C_i(\boldsymbol{\alpha}) A_i^{1/2}(\boldsymbol{\beta})$; (7) update the estimate for $\boldsymbol{\beta}$:

$$\boldsymbol{\beta}^{(t+1)} = \boldsymbol{\beta}^{(t)} - \left[ \sum_{i=1}^{N} \frac{\partial \boldsymbol{\mu}_i'}{\partial \boldsymbol{\beta}} V_i^{-1} \frac{\partial \boldsymbol{\mu}_i}{\partial \boldsymbol{\beta}} \right]^{-1} \\ \cdot \left[ \sum_{i=1}^{N} \frac{\partial \boldsymbol{\mu}_i'}{\partial \boldsymbol{\beta}} V_i^{-1} (\mathbf{y}_i - \boldsymbol{\mu}_i) \right].$$

Steps 2–7 are iterated until convergence. To illustrate step 4, consider compound symmetry, in which case the correlation is estimated by

$$\widehat{\alpha} = \frac{1}{N} \sum_{i=1}^{N} \frac{1}{n_i(n_i - 1)} \sum_{j \neq k} e_{ij} e_{ik}.$$

3.1.2 *Weighted generalized estimating equations.* As Liang and Zeger (1986) pointed out, GEE-based inferences are valid only under MCAR, due to the fact that they are based on frequentist considerations. An important exception mentioned by these authors is the situation where the working correlation structure (discussed in the previous section) happens to be correct, since then the estimates and model-based standard errors are valid under the weaker MAR. This is because then the estimating equations can be interpreted as likelihood equations. In general, of course, the working correlation structure will not be correctly specified. The ability to do so is the core motivation of the method, and therefore Robins, Rotnitzky and Zhao (1995) proposed a class of *weighted estimating equations* to allow for MAR, extending GEE.

The idea is to weight each subject's contribution in the GEEs by the inverse probability that a subject drops out at the time he dropped out. This can be calculated, for example, as

$$\nu_{id_i} \equiv P[D_i = d_i] \\ = \prod_{k=2}^{d_i - 1} (1 - P[R_{ik} = 0 | R_{i2} = \cdots = R_{i,k-1} = 1]) \\ \cdot P[R_{id_i} = 0 | R_{i2} = \cdots = R_{i,d_i-1} = 1]^{I\{d_i \leq T\}}.$$

Recall that we partitioned $\mathbf{Y}_i$ into the unobserved components $\mathbf{Y}_i^m$ and the observed components $\mathbf{Y}_i^o$. Similarly, we can make the exact same partition of $\boldsymbol{\mu}_i$ into $\boldsymbol{\mu}_i^m$ and $\boldsymbol{\mu}_i^o$. In the WGEE approach, which is proposed to reduce possible bias of $\hat{\boldsymbol{\beta}}$, the score equations to be solved when taking into account the correlation structure are

$$\begin{aligned} S(\boldsymbol{\beta}) &= \sum_{i=1}^{N} \frac{1}{\nu_{id_i}} \frac{\partial \boldsymbol{\mu}_i}{\partial \boldsymbol{\beta}'} (A_i^{1/2} C_i A_i^{1/2})^{-1} (\mathbf{y}_i - \boldsymbol{\mu}_i) \\ &= \sum_{i=1}^{N} \sum_{d=2}^{n+1} \frac{I(D_i = d)}{\nu_{id}} \\ &\quad \cdot \frac{\partial \boldsymbol{\mu}_i}{\partial \boldsymbol{\beta}'}(d) (A_i^{1/2} C_i A_i^{1/2})^{-1} \\ &\quad \cdot (d)(\mathbf{y}(d) - \boldsymbol{\mu}_i(d)) = \mathbf{0}, \end{aligned} \quad (8)$$

where $\mathbf{y}_i(d)$ and $\boldsymbol{\mu}_i(d)$ are the first $d-1$ elements of $\mathbf{y}_i$ and $\boldsymbol{\mu}_i$, respectively. We define $\partial \boldsymbol{\mu}_i / \partial \boldsymbol{\beta}'(d)$ and $(A_i^{1/2} C_i A_i^{1/2})^{-1}(d)$ analogously.

It is worthwhile to note that the recently proposed so-called *doubly robust* method (van der Laan and Robins, 2003) is more efficient and robust to a wider class of deviations. However, it is harder to implement than the original proposal. An alternative mode of analysis, generally overlooked but proposed by Schafer (2003), consists of multiply imputing the missing outcomes using a parametric model (e.g., of a random-effects or conditional type), followed by conventional GEE and conventional multiple-imputation inference on the so-completed sets of data.

### 3.2 Random-Effects Models

Unlike for correlated Gaussian outcomes, the parameters of the random-effects and population-averaged models for correlated binary data describe different types of effects of the covariates on the response probabilities (Neuhaus, 1992). Therefore, the choice between population-averaged and random-effects strategies should heavily depend on the scientific goals. Population-averaged models evaluate the success probability as a function of covariates only. With a subject-specific approach, the response is modeled as a function of covariates and parameters, specific to the subject. In such models, interpretation of fixed-effects parameters is conditional on a constant level of the random-effects parameter. Population-averaged comparisons, on the other hand, make no use of within cluster comparisons for cluster varying covariates and are therefore not useful to assess within-subject effects (Neuhaus, Kalbfleisch and Hauck, 1991). While several nonequivalent random-effects models exist, one of the most popular ones is the *generalized linear mixed model* (Breslow and Clayton, 1993), implemented in the SAS procedure NLMIXED. We will focus on this one.

3.2.1 *Generalized linear mixed models.* A general formulation of mixed-effects models is as follows. Assume that $\mathbf{Y}_i$ (possibly appropriately transformed) satisfies

$$(9) \qquad \mathbf{Y}_i|\mathbf{b}_i \sim F_i(\boldsymbol{\theta}, \mathbf{b}_i),$$

that is, conditional on $\mathbf{b}_i$, $\mathbf{Y}_i$ follows a prespecified distribution $F_i$, possibly depending on covariates and parameterized through a vector $\boldsymbol{\theta}$ of unknown parameters common to all subjects. Furthermore $\mathbf{b}_i$ is a $q$-dimensional vector of subject-specific parameters, called random effects, assumed to follow a so-called mixing distribution $G$ which may depend on a vector $\boldsymbol{\psi}$ of unknown parameters [i.e., $\mathbf{b}_i \sim G(\boldsymbol{\psi})$]. The $\mathbf{b}_i$ reflect the between-unit heterogeneity in the population with respect to the distribution of $\mathbf{Y}_i$. In the presence of random effects, conditional independence is often assumed, under which the components $Y_{ij}$ in $\mathbf{Y}_i$ are independent, conditional on $\mathbf{b}_i$. The distribution function $F_i$ in (9) then becomes a product over the $n_i$ independent elements in $\mathbf{Y}_i$.

In general, unless a fully Bayesian approach is followed, inference is based on the marginal model for $\mathbf{Y}_i$ which is obtained by integrating out the random effects over their distribution $G(\boldsymbol{\psi})$. Let $f_i(\mathbf{y}_i|\mathbf{b}_i)$ and $g(\mathbf{b}_i)$ denote the density functions that correspond to the distributions $F_i$ and $G$, respectively. We have that the marginal density function of $\mathbf{Y}_i$ equals

$$(10) \qquad f_i(\mathbf{y}_i) = \int f_i(\mathbf{y}_i|\mathbf{b}_i) g(\mathbf{b}_i) \, d\mathbf{b}_i,$$

which depends on the unknown parameters $\boldsymbol{\theta}$ and $\boldsymbol{\psi}$. Assuming independence of the units, estimates of $\widehat{\boldsymbol{\theta}}$ and $\widehat{\boldsymbol{\psi}}$ can be obtained by maximizing the likelihood function built from (10), and inferences immediately follow from classical maximum likelihood theory.

It is important to realize that the random-effects distribution $G$ is crucial in the calculation of the marginal model (10). One often assumes $G$ to be of a specific parametric form, such as a (multivariate) normal. Depending on $F_i$ and $G$, the integration in (10) may or may not be possible analytically. Proposed solutions are based on Taylor series expansions of $f_i(\mathbf{y}_i|\mathbf{b}_i)$ or on numerical approximations of the integral, such as (adaptive) Gaussian quadrature.

Note that there is an important difference with respect to the interpretation of the fixed effects $\boldsymbol{\beta}$. Under the classical linear mixed model (Verbeke and Molenberghs, 2000), we have that $E(\mathbf{Y}_i)$ equals $X_i\boldsymbol{\beta}$, such that the fixed effects have a subject-specific as well as a population-averaged interpretation. Under nonlinear mixed models, however, this no longer holds in general. The fixed effects now only reflect the conditional effect of covariates, and the marginal effect is no longer easily obtained as $E(\mathbf{Y}_i)$ is given by

$$E(\mathbf{Y}_i) = \int \mathbf{y}_i \int f_i(\mathbf{y}_i|\mathbf{b}_i) g(\mathbf{b}_i) \, d\mathbf{b}_i \, d\mathbf{y}_i.$$

However, in a biopharmaceutical context, one is often primarily interested in hypothesis testing and the random-effects framework can be used to this effect.

A general formulation of GLMM is as follows. Conditionally on random effects $\mathbf{b}_i$, it assumes that the elements $Y_{ij}$ of $\mathbf{Y}_i$ are independent, with density function usually based on a classical exponential family formulation, that is, with mean $E(Y_{ij}|\mathbf{b}_i) = a'(\eta_{ij}) = \mu_{ij}(\mathbf{b}_i)$ and variance $\text{Var}(Y_{ij}|\mathbf{b}_i) = \phi a''(\eta_{ij})$, and where, apart from a link function $h$ (e.g., the logit link for binary data or the Poisson link for counts), a linear regression model with parameters $\boldsymbol{\beta}$ and $\mathbf{b}_i$ is used for the mean [i.e., $h(\boldsymbol{\mu}_i(\mathbf{b}_i)) = X_i\boldsymbol{\beta} + Z_i\mathbf{b}_i$]. Note that the linear mixed model is a special case with identity link function. The random effects $\mathbf{b}_i$ are again assumed to be sampled from a (multivariate) normal distribution with mean $\mathbf{0}$ and covariance matrix $D$. Usually, the canonical link function is used, that is, $h = a'^{-1}$, such that $\boldsymbol{\eta}_i = X_i\boldsymbol{\beta} + Z_i\mathbf{b}_i$. When the link function is chosen to be of the logit form and the random effects are assumed to be normally distributed, the familiar logistic-linear GLMM follows.

### 3.3 Marginal versus Random-Effects Models

It is useful to underscore the difference between both model families, as well as the nature of this difference. To see the nature of the difference, consider a binary outcome variable and assume a random-intercept logistic model with linear predictor $\text{logit}[P(Y_{ij}=1|t_{ij},b_i)] = \beta_0 + b_i + \beta_1 t_{ij}$, where $t_{ij}$ is the time covariate. The conditional means $E(Y_{ij}|b_i)$, as functions of $t_{ij}$, are given by

$$(11) \qquad E(Y_{ij}|b_i) = \frac{\exp(\beta_0 + b_i + \beta_1 t_{ij})}{1 + \exp(\beta_0 + b_i + \beta_1 t_{ij})},$$

whereas the marginal average evolution is obtained by averaging over the random effects,

$$
\begin{aligned}
E(Y_{ij}) &= E[E(Y_{ij}|b_i)] \\
&= E\left[\frac{\exp(\beta_0 + b_i + \beta_1 t_{ij})}{1 + \exp(\beta_0 + b_i + \beta_1 t_{ij})}\right] \\
&\neq \frac{\exp(\beta_0 + \beta_1 t_{ij})}{1 + \exp(\beta_0 + \beta_1 t_{ij})}.
\end{aligned}
\tag{12}
$$

A graphical representation of both (11) and (12) is given in Figure 1. This implies that the interpretation of the parameters in both types of models is completely different. A schematic display is given in Figure 2. Depending on the model family (marginal or random effects), one is led to either marginal or hierarchical inference. It is important to realize that in the general case the parameter $\boldsymbol{\beta}^{\mathrm{M}}$ that results from a marginal model is different from the parameter $\boldsymbol{\beta}^{\mathrm{RE}}$ even when the latter is estimated using marginal inference. Some of the confusion surrounding this issue may result from the equality of these parameters in the very special linear mixed model case. When a random-effects model is considered, the marginal mean profile can be derived, but it will generally not produce a simple parametric form. In Figure 2 this is indicated by putting the corresponding parameter within quotes.

As an important example, consider our GLMM with logit link function, where the only random effects are intercepts $b_i$. It can then be shown that the marginal mean $\boldsymbol{\mu}_i = E(Y_{ij})$ satisfies $h(\boldsymbol{\mu}_i) \approx X_i \boldsymbol{\beta}^{\mathrm{M}}$ with

$$
\frac{\boldsymbol{\beta}^{\mathrm{RE}}}{\boldsymbol{\beta}^{\mathrm{M}}} = \sqrt{c^2 \sigma^2 + 1} > 1,
\tag{13}
$$

in which $c$ equals $16\sqrt{3}/15\pi$. Hence, although the parameters $\boldsymbol{\beta}^{\mathrm{RE}}$ in the generalized linear mixed model have no marginal interpretation, they do show a strong relationship to their marginal counterparts. Note that, as a consequence of this relationship, larger covariate effects are obtained under the random-effects model in comparison to the marginal model.

## 4. ANALYSIS OF FIRST CASE STUDY

Let us now analyze the motivating clinical trial. Therapies are recorded as A1 for primary dose of experimental drug, while B refers to nonexperimental drug and C refers to placebo. The primary contrast is between A1 and C. Emphasis is on the difference between arms at the end of the study. A graphical representation of the dropout, per study and per arm, is given in Figure 3.

The primary null hypothesis (zero difference between the treatments and placebo in terms of proportion of the $\mathrm{HAMD}_{17}$ total score above the level of 7) will be tested using both marginal models (GEE and WGEE) and random-effects models (GLMM). According to the study protocol, the models will include the fixed categorical effects of treatment, visit and treatment-by-visit interaction, as well as the continuous, fixed covariates of baseline score and baseline score-by-visit interaction. A random intercept will be included when considering the random-effects models. Analyses will be implemented using the SAS procedures GENMOD and NLMIXED.

Missing data will be handled in three different ways: (1) imputation using LOCF, (2) deletion of incomplete profiles, leading to a CC, and (3) analyzing the data as they are, consistent with ignorability (for GLMM and WGEE). A fully longitudinal approach (View 1) is considered in Section 4.1. Section 4.2 compares the results of the marginal and random-effects models. Section 4.3 focuses on Views 2 (treatment effect at last planned occasion) and 3 (last measurement obtained).

### 4.1 View 1: Longitudinal Analysis

4.1.1 *Marginal models.* First, let us consider the GEE approach. Within the SAS procedure GENMOD, the exchangeable working correlation matrix is used.

An inspection of parameter estimates and standard errors (not shown) reveals that the interaction between treatment and time is nonsignificant. The same holds in the analyses that will be done subsequently. At first sight, this suggests model simplification. However, there are a few reasons to prefer a different route. First, as stated before, a longitudinal model used in a regulatory, controlled environment is ideally sufficiently generally specified to avoid driving conclusions through models that are too simple. Sticking to a single, prespecified model also avoids dangers associated with model selection (e.g., inflated type I errors) recently reported in the literature (Hjort and Claeskens, 2003). Second, a general model allows for, as a by-product, assessment of treatment effect at the last planned occasion. Third, one can still assess the important null hypothesis of (1) no average treatment effect and (2) no treatment effect at any of the measurement occasions. These tests have been conducted and are reported in Table 2; for (1), also the estimated average treatment effect is reported.

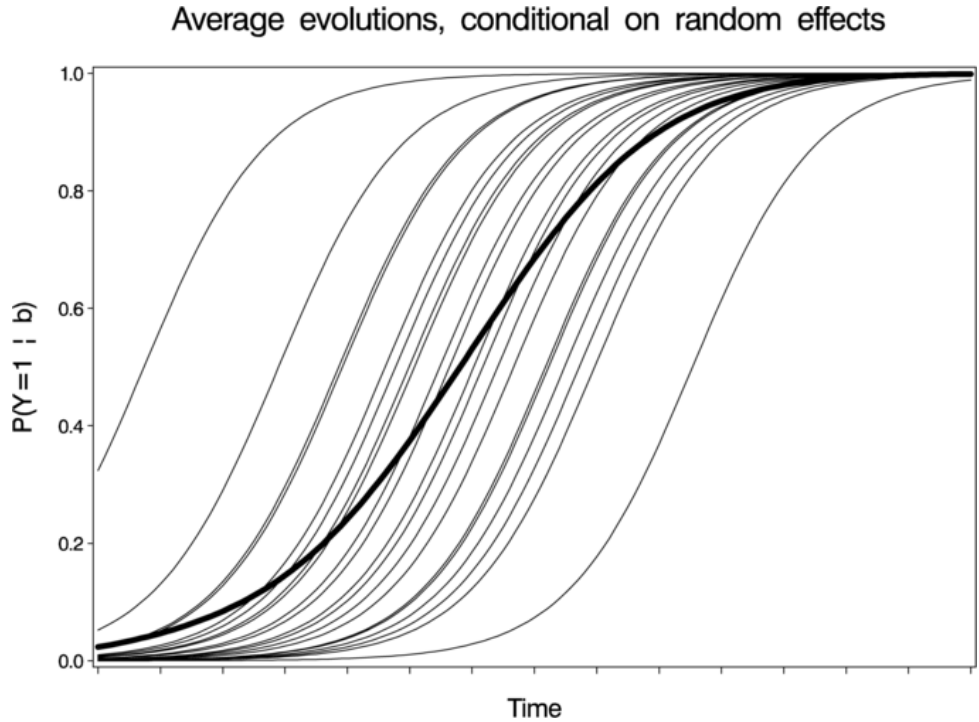

Fig. 1. *Graphical representation of a random-intercept logistic curve across a range of levels of the random intercept, together with the corresponding marginal curve.*

In many cases, the empirically corrected standard errors are larger than the model-based ones. This is because model-based standard errors are the ones that would be obtained if the estimating equations were true likelihood equations, that is, when the working correlation structure is correct. In such cases likelihood inference enjoys optimality. However, since the working correlation structure is allowed to be misspecified, model-based standard errors will be bi-ased and it is advisable to base conclusions on empirically corrected standard errors.

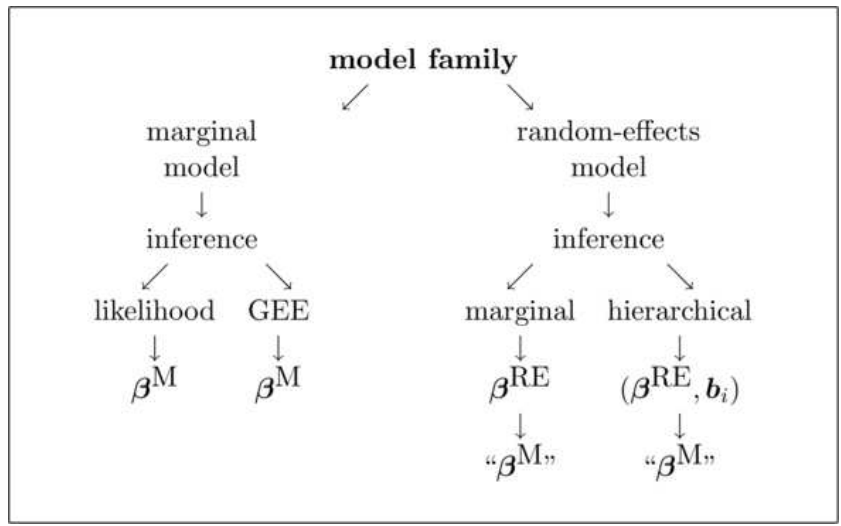

Fig. 2. *Representation of model families and corresponding inference. A superscript* M *stands for marginal;* RE *denotes random effects. A parameter within quotes indicates that marginal functions but no direct marginal parameters are obtained.*

Turning to WGEE, the method is applied to perform an analysis that is correct under MAR, not only under MCAR as in ordinary GEE. This procedure is a bit more involved in terms of fitting the model to the data. We will outline the main steps. The SAS code is available from the authors upon request.

To compute the necessary weights, we first fit the dropout model using logistic regression. The outcome `drop` is binary and indicates whether or not dropout occurs at a given time. The response value at the previous occasion (`prevhamd`) and treatment are included as covariates. Next, the predicted probabilities of dropout are translated into weights, defined at the individual measurement level. Let us describe the procedure to construct the inverse weights. At the first occasion define $w_{i1} = 1$. At other than the last occasion, the quantity of interest equals the cumulative weight over the previous occasions, multiplied by (1 − the predicted probability of dropout). At the last occasion *within a sequence where dropout occurs*, it is multiplied by the predicted probability of dropout. At the end of the process this quantity is inverted to yield the actual weight. After these preparations we merely need to include the weights by means of the `scwgt` statement within the GENMOD procedure. Together with the use of the `repeated` statement, WGEE follows. Also here we use the exchangeable working correlation matrix.

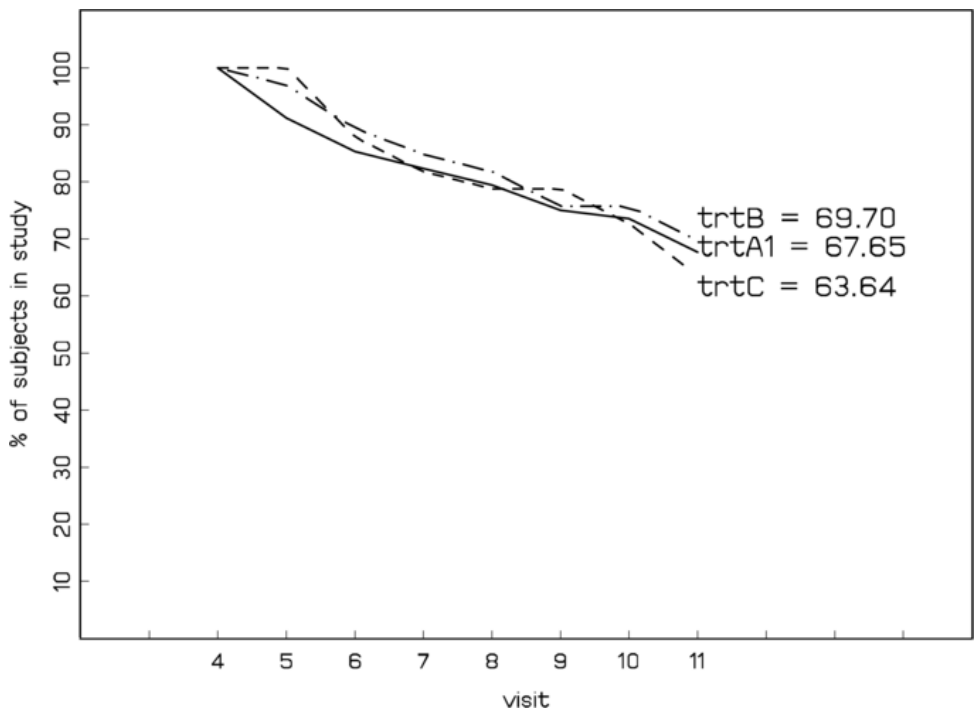

FIG. 3. *Evolution of dropout per study and per treatment arm. Treatment arms A1 and C, being the ones of primary interest, are shown by the solid and dashed lines, respectively.*

TABLE 2
*Depression trial, View 1: GEE, WGEE and GLMM. Tests for* (1) *the joint null hypothesis of no treatment effect at none of the time points and* (2) *the hypothesis of no average treatment effect*

|  | Joint effects | | | Mean effects | | |
| --- | --- | --- | --- | --- | --- | --- |
|  | A1<br>(8 d.f.)<br>$p$ | B<br>(8 d.f.)<br>$p$ | A1 & B1<br>(16 d.f.)<br>$p$ | A1<br>(1 d.f.)<br>$p$ (est., s.e.) | B<br>(1 d.f.)<br>$p$ (est., s.e.) | A1 & B1<br>(2 d.f.)<br>$p$ |
| Analysis |  |  |  |  |  |  |
| CC (GEE) | 0.4278 | 0.9859 | 0.8444 | 0.5259 (−2.66; 4.19) | 0.9165 (0.47; 4.50) | 0.5845 |
| LOCF (GEE) | 0.7008 | 0.9956 | 0.9768 | 0.7713 (−1.15; 3.96) | 0.9070 (0.49; 4.19) | 0.8605 |
| MAR (GEE) | 0.6465 | 0.9931 | 0.9413 | 0.6015 (−1.92; 3.67) | 0.8671 (0.65; 3.89) | 0.6804 |
| MAR (WGEE) | 0.1690 | 0.7601 | 0.5372 | 0.5477 (2.61; 4.34) | 0.3883 (3.97; 4.60) | 0.7224 |
| CC (GLMM) | 0.7572 | 0.9743 | 0.7233 | 0.4954 (−0.40; 0.59) | 0.2671 (0.64; 0.57) | 0.0440 |
| LOCF (GLMM) | 0.7363 | 0.9953 | 0.9763 | 0.1571 (−0.66; 0.47) | 0.4555 (−0.34; 0.45) | 0.3611 |
| MAR (GLMM) | 0.7476 | 0.9738 | 0.7152 | 0.4495 (−0.41; 0.55) | 0.2844 (0.58; 0.54) | 0.0375 |

Let us now turn to the results. The marginal models reveal nonsignificant treatment effects in all cases, for either the composite hypothesis of no treatment effects or the hypotheses of no average effects. This holds for both arms separately, as well as for the two arms jointly. Corresponding to the 1-degree-of-freedom (d.f.) tests, parameter estimates and standard errors can be estimated as well. For conciseness, only empirically corrected standard errors are shown. A strong difference is observed between the WGEE and other cases. Since this is the only one valid under MAR, it is clear that there are dangers associated with methods that are too simple. Furthermore, some of the CC $p$-values are smaller than their MAR and LOCF counterparts.

4.1.2 *Random-effect models.* To fit generalized linear mixed models, we use the SAS procedure NLMIXED, which allows fitting a wide class of linear, generalized linear and nonlinear mixed models. It relies on numerical integration. Not only are different integral approximations available, the principal ones being (adaptive) Gaussian quadrature, it also includes a number of optimization algorithms. The difference between nonadaptive and adaptive Gaussian quadrature is that for the first procedure the quadrature points are centered at zero for each of the random effects and the current random-effects covariance matrix is used as the scale matrix, while for the latter the quadrature points will be appropriately centered and scaled, such that more quadrature points lie in the region of interest (Molenberghs and Verbeke, 2005). We will use both adaptive and nonadaptive quadrature, with several choices for the number of quadrature points, to check the stability of the results over a variety of choices for these numerical choices.

Precisely, we initiate the model fitting using nonadaptive Gaussian quadrature, together with the quasi-Newton optimization algorithm. The number of quadrature points is left to be determined by the procedure, and all starting values are set equal to 0.5. Using the resulting parameter estimates, we keep these choices but hold the number of quadrature points fixed (2, 3, 5, 10, 20 and 50). Subsequently, we switch to adaptive Gaussian quadrature (step 2). Finally, the quasi-Newton optimization is replaced by the Newton–Raphson optimization (step 3). The effect of the method and the number of quadrature points is graphically represented in Figure 4 for a selected parameter (treatment effect of A1). While the differences between these choices are purely numerical, we do notice differences between the results, illustrating that a numerical sensitivity analysis matters. The parameter estimates tend to stabilize with increasing number of quadrature points. However, nonadaptive Gaussian quadrature needs obviously more quadrature points than adaptive Gaussian quadrature.

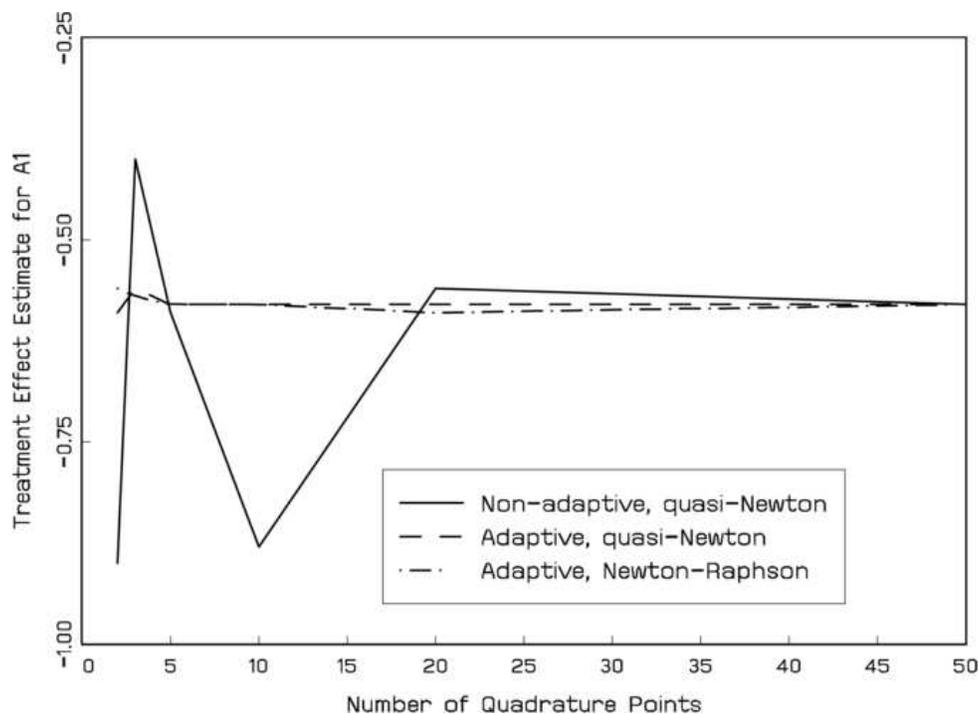

Fig. 4. *The effect of adaptive versus nonadaptive quadrature, quasi-Newton versus Newton–Raphson and the number of quadrature points on the treatment effect parameter for arm* A1.

Focusing on the results for 50 quadrature points, we observe that the parameter estimates for steps 1 and 2 are the same. On the other hand, parameter estimates for step 3 are different (order of $10^{-3}$, visible in $p$-values). In spite of the differences in parameter estimates, it is noteworthy that the likelihood is the same in all steps, due to a flat likelihood. This was confirmed by running all steps again using the parameter estimates of step 3 as starting values, at which point the parameter estimates all coincide. Thus, it may happen that the optimization routine has only seemingly converged.

Exactly as in the marginal model case, we assessed average treatment effect as well as treatment effect at any of the times. The results are reported in Table 2 as well. The parameter $p$-values are more varied across methods than in the marginal model case. The most striking feature is that there is evidence for a treatment effect in the two groups together, under MAR, and also with CC. Note that the corresponding 1-degree-of-freedom tests do not show significance. In the LOCF case some $p$-values are smaller, while others are larger. This contradicts a common belief that LOCF is conservative. Molenberghs et al. (2004) have shown that both conservative and liberal behavior is possible.

### 4.2 Marginal versus Random-Effects Models

In all cases, the variability of the random effect (standard deviation parameter $\sigma$) is highly significant. This implies that the GEE parameters and the random-effects parameters cannot be compared directly. If the conversion factor (13) is computed, then one roughly finds a factor of about 2.5. We note that this factor is not reproduced when the two sets of estimates (estimates not shown) are directly compared. This is due to the fact that (13) operates at the true population parameter level, while we only have parameter estimates at our disposal. Since many of the estimates are only marginally or not significant, it is not unexpected, therefore, to observe deviations from this relationship, even though the general tendency is preserved in most cases.

### 4.3 Views 2 and 3: Single Time Point Analysis

When emphasis is on the last measurement occasion, LOCF and CC are straightforward to use. When the last observed measurement is of interest (a different scientific question), the analysis is not different from the one obtained under LOCF, but, of course, in this case CC is not an option.

Since the outcome is a dichotomous response, the data can be summarized in a $2 \times k$ table, where $k$ represents the number of treatments. The analysis essentially consists of comparing the proportions of success or failure in all groups. For this purpose, both Pearson's chi-squared test (Agresti, 1990) and Fisher's exact test (Freeman and Halton, 1951) will be used. Nevertheless, it is still possible to obtain inferences from a full longitudinal model in this context. We add these for the sake of reference, but it should be understood that the analysis using a simple model for the last time point only is more in line with practice. When an ignorable analysis is considered, one has to explicitly consider all incomplete profiles in order to correctly incorporate all information available. Thus, one has to consider a longitudinal model.

Placebo C is considered as the reference treatment. Let $\alpha_i$ be the effect of treatment arm $i$ at the last measurement occasion, where $i =$ A1, B or C. We wish to test whether, at the last measurement occasion, all treatment effects are equal. This translates into $\alpha_{A1} = \alpha_B = \alpha_C$ or, equivalently, into $\alpha_{A1} - \alpha_C = \alpha_B - \alpha_C = 0$. Such contrasts can be obtained very easily using the SAS procedure NLMIXED. Table 3 shows a summary of the results in terms of $p$-values.

The GLMMs lead to a small difference between CC and MAR: both are borderline. On the other hand, the GLMM for LOCF leads to a nonsignificant result. An endpoint analysis (i.e., using the last available measurement) shows the same result for LOCF (nonsignificant), whereas the result for CC becomes significant. An endpoint analysis leads to a completely different picture, with results that are strongly different (significant) from the GLMM model. This illustrates that the choice between modeling techniques is far from an academic question, but can have profound impact on the study conclusions, ranging from highly significant over borderline (non)significant to highly nonsignificant.

TABLE 4
*Age related macular degeneration trial. Parameter estimates (model-based standard errors; empirically corrected standard errors) for the marginal models: GEE on the CC and LOCF population, and on the observed data; in the latter case WGEE is also used*

| Effect | Parameter | CC | LOCF | Observed data Unweighted | WGEE |
|---|---|---|---|---|---|
| Int. 4 | $\beta_{11}$ | $-1.01\ (0.24; 0.24)$ | $-0.87\ (0.20; 0.21)$ | $-0.87\ (0.21; 0.21)$ | $-0.98\ (0.10; 0.44)$ |
| Int. 12 | $\beta_{21}$ | $-0.89\ (0.24; 0.24)$ | $-0.97\ (0.21; 0.21)$ | $-1.01\ (0.21; 0.21)$ | $-1.78\ (0.15; 0.38)$ |
| Int. 24 | $\beta_{31}$ | $-1.13\ (0.25; 0.25)$ | $-1.05\ (0.21; 0.21)$ | $-1.07\ (0.22; 0.22)$ | $-1.11\ (0.15; 0.33)$ |
| Int. 52 | $\beta_{41}$ | $-1.64\ (0.29; 0.29)$ | $-1.51\ (0.24; 0.24)$ | $-1.71\ (0.29; 0.29)$ | $-1.72\ (0.25; 0.39)$ |
| Trt. 4 | $\beta_{12}$ | $0.40\ (0.32; 0.32)$ | $0.22\ (0.28; 0.28)$ | $0.22\ (0.28; 0.28)$ | $0.80\ (0.15; 0.67)$ |
| Trt. 12 | $\beta_{22}$ | $0.49\ (0.31; 0.31)$ | $0.55\ (0.28; 0.28)$ | $0.61\ (0.29; 0.29)$ | $1.87\ (0.19; 0.61)$ |
| Trt. 24 | $\beta_{32}$ | $0.48\ (0.33; 0.33)$ | $0.42\ (0.29; 0.29)$ | $0.44\ (0.30; 0.30)$ | $0.73\ (0.20; 0.52)$ |
| Trt. 52 | $\beta_{42}$ | $0.40\ (0.38; 0.38)$ | $0.34\ (0.32; 0.32)$ | $0.44\ (0.37; 0.37)$ | $0.74\ (0.31; 0.52)$ |
| Corr. | $\rho$ | 0.39 | 0.44 | 0.39 | 0.33 |

TABLE 3
*Depression trial, Views 2 and 3: p-values are reported (mixed refers to the assessment of treatment at the last visit based on a generalized linear mixed model)*

| Method | Model | p-Value |
|---|---|---|
| CC | Mixed | 0.0463 |
|  | Pearson's chi-squared test | 0.0357 |
|  | Fisher's exact test | 0.0336 |
| LOCF | Mixed | 0.1393 |
|  | Pearson's chi-squared test | 0.1553 |
|  | Fisher's exact test | 0.1553 |
| MAR | Mixed | 0.0500 |

## 5. ANALYSIS OF THE SECOND CASE STUDY

We compare analyses performed on the completers only (CC), on the LOCF imputed data and on the observed data. For the observed, partially incomplete data, GEE is supplemented with WGEE. Furthermore, a random-intercepts GLMM is considered, based on numerical integration. The GEE analyses are reported in Table 4 and the random-effects models in Table 5. For GEE, a working exchangeable correlation matrix is considered. The model has four intercepts and four treatment effects. To be precise, the marginal regression model takes the form

$$\text{logit}[P(Y_{ij}=1|T_i)] = \beta_{j1} + \beta_{j2}T_i,$$

where $j=1,\ldots,4$ refers to measurement occasion, $T_i$ is the treatment assignment for subject $i=1,\ldots,240$ and $Y_{ij}$ is the indicator for whether or not three lines of vision have been lost for subject $i$ at time $j$. The advantage of having separate treatment effects at each time is that particular attention can be given to the treatment effect assessment at the last planned measurement occasion (i.e., after one year). From Table 4 it is clear that the model-based and empirically corrected standard errors agree extremely well. This is due to the unstructured nature of the full time by treatment mean structure. However, we do observe differences in the WGEE analyses. Not only are the parameter estimates mildly different between the two GEE versions, but there is a dramatic difference between the model-based and empirically corrected standard errors. Nevertheless, the two sets of empirically corrected standard errors agree very closely, which is reassuring.

TABLE 5
*Age related macular degeneration trial. Parameter estimates (standard errors) for the random-intercept models: Numerical-integration-based fits (adaptive Gaussian quadrature) on the CC and LOCF population, and on the observed data (direct-likelihood)*

| Effect | Parameter | CC | LOCF | Direct likelihood |
|---|---|---|---|---|
| Int. 4 | $\beta_{11}$ | $-1.73\ (0.42)$ | $-1.63\ (0.39)$ | $-1.50\ (0.36)$ |
| Int. 12 | $\beta_{21}$ | $-1.53\ (0.41)$ | $-1.80\ (0.39)$ | $-1.73\ (0.37)$ |
| Int. 24 | $\beta_{31}$ | $-1.93\ (0.43)$ | $-1.96\ (0.40)$ | $-1.83\ (0.39)$ |
| Int. 52 | $\beta_{41}$ | $-2.74\ (0.48)$ | $-2.76\ (0.44)$ | $-2.85\ (0.47)$ |
| Trt. 4 | $\beta_{12}$ | $0.64\ (0.54)$ | $0.38\ (0.52)$ | $0.34\ (0.48)$ |
| Trt. 12 | $\beta_{22}$ | $0.81\ (0.53)$ | $0.98\ (0.52)$ | $1.00\ (0.49)$ |
| Trt. 24 | $\beta_{32}$ | $0.77\ (0.55)$ | $0.74\ (0.52)$ | $0.69\ (0.50)$ |
| Trt. 52 | $\beta_{42}$ | $0.60\ (0.59)$ | $0.57\ (0.56)$ | $0.64\ (0.58)$ |
| R.I. s.d. | $\tau$ | $2.19\ (0.27)$ | $2.47\ (0.27)$ | $2.20\ (0.25)$ |
| R.I. var. | $\tau^2$ | $4.80\ (1.17)$ | $6.08\ (1.32)$ | $4.83\ (1.11)$ |

When comparing parameter estimates across CC, LOCF and observed data analyses, it is clear that LOCF has the effect of artificially increasing the correlation between measurements. The effect is mild in this case. The parameter estimates of the observed-data GEE are close to the LOCF results for earlier time points and close to CC for later time points. This is to be expected, because at the start of the study the LOCF and observed populations are virtually the same, with the same holding between CC and observed populations near the end of the study. Note also that the treatment effect under LOCF, especially at 12 weeks and after 1 year, is biased downward in comparison to the GEE analyses. To properly use the information in the missingness process, WGEE can be used. To this end, a logistic regression for dropout, given covariates and previous outcomes, needs to be fitted. Parameter estimates and standard errors are given in Table 6. Intermittent missingness will be ignored. Covariates of importance are treatment assignment, the level of lesions at baseline (a four-point categorical variable, for which three dummies are needed) and time at which dropout occurs. For the latter covariates, there are three levels, since dropout can occur at times 2, 3 or 4. Hence, two dummy variables are included. Finally, the previous outcome does not have a significant impact, but will be kept in the model nevertheless. In spite of there being no strong evidence for MAR, the results between GEE and WGEE differ quite a bit. It is noteworthy that at 12 weeks, a treatment effect is observed with WGEE that goes unnoticed with the other marginal analyses. This finding is mildly confirmed by the random-intercept model when the data as observed are used.

The results for the random-effects models are given in Table 5. We observe the usual relationship between the marginal parameters of Table 4 and their random-effects counterparts. Note also that the random-intercepts variance is largest under LOCF, underscoring again that this method artificially increases the association between measurements on the same subject. In this case, unlike for the marginal models, LOCF and in fact also CC slightly to considerably overestimate the treatment effect at certain times, in particular at 4 and 24 weeks.

TABLE 6
*Age related macular degeneration trial. Parameter estimates (standard errors) for a logistic regression model to describe dropout*

| Effect | Parameter | Estimate (s.e.) |
|---|---|---|
| Intercept | $\psi_0$ | 0.14 (0.49) |
| Previous outcome | $\psi_1$ | 0.04 (0.38) |
| Treatment | $\psi_2$ | −0.86 (0.37) |
| Lesion level 1 | $\psi_{31}$ | −1.85 (0.49) |
| Lesion level 2 | $\psi_{32}$ | −1.91 (0.52) |
| Lesion level 3 | $\psi_{33}$ | −2.80 (0.72) |
| Time 2 | $\psi_{41}$ | −1.75 (0.49) |
| Time 3 | $\psi_{42}$ | −1.38 (0.44) |

## 6. DISCUSSION

In this paper we have indicated that a variety of approaches is possible when analyzing incomplete longitudinal data from clinical trials. First, unlike in the continuous case where the linear mixed model is the main mode of analysis, one has the choice between a marginal model (generalized estimating equations, GEE) and a random-effects approach (generalized linear mixed models, GLMM). While these may provide similar results in terms of hypothesis testing, things are different when the models are used for estimation purposes, because the parameters have quite different meanings. Both GEE and GLMM can be used when data are incomplete. For GLMM this holds under the fairly general assumption of an MAR mechanism, while for GEE the stronger MCAR is required. However, GEE can be extended to weighted GEE, making it also valid under MAR. Current statistical computing power has brought both GLMM and WGEE within reach, and we have implemented such analyses in the real-life setting of clinical trials on depression and on macular degeneration. This underscores that simple but potentially highly restrictive modes of analyses, such as CC or LOCF, should no longer be seen as the preferred mode of analysis. This message is in line with the one reached for continuous outcomes (Molenberghs et al., 2004).

While in the studies considered here there are no extreme differences between the various analyses conducted, some differences are noticeable, especially in the second case study (Molenberghs et al., 2004). So, generally, caution should be used and it is best to move away from the overly simple methods.

Note that such a full longitudinal approach under MAR is also very sensible, even when one is interested in an effect at one particular scheduled measurement occasion, say, the treatment effect at the last scheduled visit. Indeed, an ignorable analysis takes all information into account, not only from complete observations, but also from incomplete ones, through the conditional expectation of the missing measurements given the observed ones. Thus, when combined with an analysis where the treatment allocation is used "as randomized" rather than "as treated," such an approach is fully compatible with the intention-to-treat principle.

When there is residual doubt about the plausibility of MAR, one can conduct a sensitivity analysis. Many proposals have been made, but this remains an active area of research. Obviously, a number of MNAR models can be fitted, provided one is prepared to approach formal aspects of model comparison with due caution. Such analyses can be complemented with appropriate (global and/or local) influence analyses. Some sensitivity analyses frameworks have been provided by Rotnitzky, Robins and Scharfstein (1998), Forster and Smith (1998), Raab and Donnelly (1999), Kenward, Goetghebeur and Molenberghs (2001), van Steen, Molenberghs, Verbeke and Thijs (2003) and Jansen et al. (2003).

## ACKNOWLEDGMENTS

Ivy Jansen, Caroline Beunckens and Geert Molenberghs gratefully acknowledge support from Fonds Wetenschappelijk Onderzoek-Vlaanderen Research Project G.0002.98 (Sensitivity Analysis for Incomplete and Coarse Data) and from the Belgian IUAP/PAI network (Statistical Techniques and Modeling for Complex Substantive Questions with Complex Data). We thank Eli Lilly for kind permission to use their data.

## REFERENCES

Aerts, M., Geys, H., Molenberghs, G. and Ryan, L. M. (2002). *Topics in Modelling of Clustered Data*. Chapman and Hall, London. MR2022882

Afifi, A. and Elashoff, R. (1966). Missing observations in multivariate statistics. I. Review of the literature. *J. Amer. Statist. Assoc.* **61** 595–604. MR0203865

Agresti, A. (1990). *Categorical Data Analysis*. Wiley, New York. MR1044993

Ashford, J. R. and Sowden, R. R. (1970). Multivariate probit analysis. *Biometrics* **26** 535–546.

Bahadur, R. R. (1961). A representation of the joint distribution of responses to $n$ dichotomous items. In *Studies in Item Analysis and Prediction* (H. Solomon, ed.) 169–176. Stanford Univ. Press, Stanford, CA. MR0121894

Bowman, D. and George, E. O. (1995). A saturated model for analyzing exchangeable binary data: Applications to clinical and developmental toxicity studies. *J. Amer. Statist. Assoc.* **90** 871–879.

Breslow, N. E. and Clayton, D. G. (1993). Approximate inference in generalized linear mixed models. *J. Amer. Statist. Assoc.* **88** 9–25.

Dale, J. R. (1986). Global cross-ratio models for bivariate, discrete, ordered responses. *Biometrics* **42** 909–917.

Dempster, A. P., Laird, N. M. and Rubin, D. B. (1977). Maximum likelihood from incomplete data via the EM algorithm (with discussion). *J. Roy. Statist. Soc. Ser. B* **39** 1–38. MR0501537

Dempster, A. P. and Rubin, D. B. (1983). Overview. In *Incomplete Data in Sample Surveys* **2**. *Theory and Bibliographies* (W. G. Madow, I. Olkin and D. B. Rubin, eds.) 3–10. Academic Press, New York.

Diggle, P. J., Heagerty, P. J., Liang, K.-Y. and Zeger, S. L. (2002). *Analysis of Longitudinal Data*, 2nd ed. Oxford Univ. Press, New York. MR2049007

Ekholm, A. (1991). Algorithms versus models for analyzing data that contain misclassification errors (with response). *Biometrics* **47** 1171–1182.

Fahrmeir, L. and Tutz, G. (2001). *Multivariate Statistical Modelling Based on Generalized Linear Models*. Springer, Heidelberg. MR1832899

Fitzmaurice, G. M. (2003). Methods for handling dropouts in longitudinal clinical trials. *Statist. Neerlandica* **57** 75–99. MR2055522

Forster, J. J. and Smith, P. W. F. (1998). Model-based inference for categorical survey data subject to non-ignorable non-response. *J. R. Stat. Soc. Ser. B Stat. Methodol.* **60** 57–70. MR1625664

Freeman, G. H. and Halton, J. H. (1951). Note on an exact treatment of contingency, goodness of fit and other problems of significance. *Biometrika* **38** 141–149. MR0042666

Geys, H., Molenberghs, M. and Lipsitz, S. R. (1998). A note on the comparison of pseudo-likelihood and generalized estimating equations for marginal odds ratio models with exchangeable association structure. *J. Statist. Comput. Simulation* **62** 45–72.

Gilula, Z. and Haberman, S. J. (1994). Conditional log-linear models for analyzing categorical panel data. *J. Amer. Statist. Assoc.* **89** 645–656. MR1294088

Glonek, G. F. V. and McCullagh, P. (1995). Multivariate logistic models. *J. Roy. Statist. Soc. Ser. B* **55** 533–546.

Hartley, H. O. and Hocking, R. R. (1971). The analysis of incomplete data. *Biometrics* **27** 783–823.

Hjort, N. L. and Claeskens, G. (2003). Frequentist model average estimators. *J. Amer. Statist. Assoc.* **98** 879–899. MR2041481

Jansen, I., Molenberghs, G., Aerts, M., Thijs, H. and van Steen, K. (2003). A local influence approach applied

to binary data from a psychiatric study. *Biometrics* **59** 410–418.MR1987408

Kenward, M. G., Goetghebeur, E. J. T. and Molenberghs, G. (2001). Sensitivity analysis for incomplete categorical data. *Statistical Modelling* **1** 31–48.

Kenward, M. G. and Molenberghs, G. (1998). Likelihood based frequentist inference when data are missing at random. *Statist. Sci.* **13** 236–247.MR1665713

Lang, J. B. and Agresti, A. (1994). Simultaneously modeling joint and marginal distributions of multivariate categorical responses. *J. Amer. Statist. Assoc.* **89** 625–632.

le Cessie, S. and van Houwelingen, J. C. (1994). Logistic regression for correlated binary data. *Appl. Statist.* **43** 95–108.

Liang, K.-Y. and Zeger, S. L. (1986). Longitudinal data analysis using generalized linear models. *Biometrika* **73** 13–22.MR0836430

Little, R. J. A. and Rubin, D. B. (2002). *Statistical Analysis with Missing Data*, 2nd ed. Wiley, New York.MR1925014

Longford, N. (1993). Inference about variation in clustered binary data. Paper presented at the Multilevel Conference, Rand Corporation, Los Angeles.

Mallinckrodt, C. H., Clark, W. S., Carroll, R. J. and Molenberghs, G. (2003a). Assessing response profiles from incomplete longitudinal clinical trial data under regulatory considerations. *J. Biopharmaceutical Statistics* **13** 179–190.

Mallinckrodt, C. H., Clark, W. S. and David, S. R. (2001a). Type I error rates from mixed-effects model repeated measures versus fixed effects analysis of variance with missing values imputed via last observation carried forward. *Drug Information J.* **35** 1215–1225.

Mallinckrodt, C. H., Clark, W. S. and David, S. R. (2001b). Accounting for dropout bias using mixed-effects models. *J. Biopharmaceutical Statistics* **11** 9–21.

Mallinckrodt, C. H., Sanger, T. M., Dubé, S., DeBrota, D. J., Molenberghs, G., Carroll, R. J., Potter, W. Z. and Tollefson, G. D. (2003b). Assessing and interpreting treatment effects in longitudinal clinical trials with missing data. *Biological Psychiatry* **53** 754–760.

McCullagh, P. and Nelder, J. A. (1989). *Generalized Linear Models*, 2nd ed. Chapman and Hall, London.MR0727836

Molenberghs, G. and Lesaffre, E. (1994). Marginal modeling of correlated ordinal data using a multivariate Plackett distribution. *J. Amer. Statist. Assoc.* **89** 633–644.

Molenberghs, G. and Lesaffre, E. (1999). Marginal modelling of multivariate categorical data. *Statistics in Medicine* **18** 2237–2255.

Molenberghs, G., Thijs, H., Jansen, I., Beunckens, C., Kenward, M. G., Mallinckrodt, C. and Carroll, R. J. (2004). Analyzing incomplete longitudinal clinical trial data. *Biostatistics* **5** 445–464.

Molenberghs, G. and Verbeke, G. (2005). *Models for Discrete Longitudinal Data*. Springer, New York.MR2171048

Neuhaus, J. M. (1992). Statistical methods for longitudinal and clustered designs with binary responses. *Statistical Methods in Medical Research* **1** 249–273.

Neuhaus, J. M., Kalbfleisch, J. D. and Hauck, W. W. (1991). A comparison of cluster-specific and population-averaged approaches for analyzing correlated binary data. *Internat. Statist. Rev.* **59** 25–35.

Pharmacological Therapy for Macular Degeneration Study Group (1997). Interferon $\alpha$-IIA is ineffective for patients with choroidal neovascularization secondary to age-related macular degeneration. Results of a prospective randomized placebo-controlled clinical trial. *Archives of Ophthalmology* **115** 865–872.

Plackett, R. L. (1965). A class of bivariate distributions. *J. Amer. Statist. Assoc.* **60** 516–522.MR0183042

Raab, G. M. and Donnelly, C. A. (1999). Information on sexual behaviour when some data are missing. *Appl. Statist.* **48** 117–133.

Robins, J. M., Rotnitzky, A. and Zhao, L. P. (1995). Analysis of semiparametric regression models for repeated outcomes in the presence of missing data. *J. Amer. Statist. Assoc.* **90** 106–121.MR1325118

Rotnitzky, A., Robins, J. M., and Scharfstein, D. O. (1998). Semiparametric regression for repeated outcomes with nonignorable nonresponse. *J. Amer. Statist. Assoc.* **93** 1321–1339.MR1666631

Rubin, D. B. (1976). Inference and missing data. *Biometrika* **63** 581–592.MR0455196

Rubin, D. B. (1987). *Multiple Imputation for Nonresponse in Surveys*. Wiley, New York.
MR0899519

Rubin, D. B., Stern, H. S. and Vehovar, V. (1995). Handling "don't know" survey responses: The case of the Slovenian plebiscite. *J. Amer. Statist. Assoc.* **90** 822–828.

Schafer, J. (2003). Multiple imputation in multivariate problems when the imputation and analysis models differ. *Statist. Neerlandica* **57** 19–35.MR2055519

Siddiqui, O. and Ali, M. W. (1998). A comparison of the random-effects pattern mixture model with last-observation-carried-forward (LOCF) analysis in longitudinal clinical trials with dropouts. *J. Biopharmaceutical Statistics* **8** 545–563.

Stiratelli, R., Laird, N. and Ware, J. H. (1984). Random effects models for serial observations with binary response. *Biometrics* **40** 961–971.

van der Laan, M. J. and Robins, J. M. (2003). *Unified Methods for Censored Longitudinal Data and Causality*. Springer, New York.MR1958123

van Steen, K., Molenberghs, G., Verbeke, G. and Thijs, H. (2001). A local influence approach to sensitivity analysis of incomplete longitudinal ordinal data. *Statistical Modelling* **1** 125–142.

Verbeke, G. and Molenberghs, G. (2000). *Linear Mixed Models for Longitudinal Data*. Springer, New York.MR1880596

Wolfinger, R. and O'Connell, M. (1993). Generalized linear mixed models: A pseudo-likelihood approach. *J. Statist. Comput. Simulation* **48** 233–243.